\documentclass[11pt,letterpaper]{amsart}

\usepackage{amssymb}
\usepackage{amsthm}
\usepackage{amsfonts}
\usepackage{amsmath}
\usepackage{amstext}
\usepackage{graphicx}

\numberwithin{equation}{section}

\theoremstyle{plain}

\theoremstyle{definition}

\theoremstyle{remark}

\newcommand{\R}{\mathbb{R}}

\newcommand{\Z}{\mathbb{Z}}

\DeclareMathOperator{\hlink}{\tilde\Lambda }
\DeclareMathOperator{\link}{ \Lambda }

\DeclareMathOperator{\ahlink} {\tilde\Omega }

\def\({(\!(}
\def\){)\!)}

\def\accentclass@{7}
\def\makeacc@#1#2{\def#1{\mathaccent"\accentclass@#2 }}
\makeacc@\cir{017}

\date{June 1, 2001}

\title[How to show a set is not algebraic]
{How to show a set is not algebraic}

\author[C. McCrory] {Clint McCrory}
\address{Mathematics Department, University of Georgia,
Athens, Georgia 30602-7403, USA}
\email{clint@math.uga.edu}


\newcommand{\abstracttext}{We revisit Akbulut and King's first example of a
compact
semialgebraic set which satisfies Sullivan's local Euler
characteristic condition, but which is not homeomorphic to an algebraic set. A
nontrivial obstruction is computed using the link operator on the ring of
constructible functions.}

\begin{document}
\begin{abstract} \abstracttext \end{abstract}
\maketitle

There are many local topological conditions satisfied by real algebraic
sets. The
simplest of these was discovered by Sullivan  \cite{sullivan} more than thirty years ago:
The link of
every point has even Euler characteristic. In low dimensions this is the only
obstruction for a set to be algebraic. More precisely, if a compact
semialgebraic
set of dimension less than three satisfies Sullivan's condition, then it is
homeomorphic to an algebraic set. Akbulut and King
\cite{akbulut-king} found an example of a compact 3-dimensional
semialgebraic set
satisfying Sullivan's condition which is not homeomorphic to an algebraic
set, and
this is the example we will discuss.

The method we use to compute the obstruction is
due to Parusi\'nski and the author
\cite{mccrory-parusinski1}.
We have also found  an enormous list of independent obstructions
in dimension four \cite{mccrory-parusinski2}. The yoga of algebraically constructible
functions used
to prove our obstructions vanish for algebraic sets is presented in the paper of
Isabelle Bonnard in this volume.

Let $X$ be a semialgebraic set in $\R^n$,
and let $p\in X$. The {\em link} $L_p$ of
$X$ at $p$ is the intersection of $X$ with a sphere of small
radius $\epsilon>0$ centered at $p$. For sufficiently small $\epsilon$ the
topological type of $L_p$ is independent of $\epsilon$, and a
neighborhood of $p$ in $X$ is homeomorphic to the cone on $L_p$. (The link
$L_p$ is
homeomorphic to the boundary of the simplicial star of $p$ in $X$ for a
semialgebraic
triangulation of $X$.)
Sullivan's theorem  is that if $X$ is algebraic then for all
$p\in X$
the Euler characteristic $\chi(L_p)$ is even.

For example, consider the {\em Cartan umbrella} $X=\{(x,y,z)\ |\ x^2 =
zy^2\}$ in
$\R^3$ (Figure 1). For each real number $t$, let $H_t$ be the horizontal plane
$z=t$. For $t>0$, $X\cap H_t$ is the two lines $x=\pm(\sqrt t) y$, $z=t$. For
$t=0$,
$X\cap H_t$ is the line $x=0$, $z=0$. For $t<0$, $X\cap H_t$ is the point
$x=0$,
$y=0$, $z=t$.

The space $X$ has five  strata:
\begin{itemize}
\item two 2-strata $f=\{x^2=zy^2, y<0\}$, $f'=\{x^2=zy^2, y>0\}$,
\item two 1-strata $d=\{x=0, y=0, z<0\}$, $e=\{x=0, y=0, z>0\}$,
\item one 0-stratum $a =\{(0,0,0)\}$.
\end{itemize}

The local topology of $X$ is constant along each of these strata. If $p \in
X$, then
\begin{itemize}
\item if $p\in f\cup f'$, $L_p$ is topologically a circle, so $\chi(L_p)=0$;
\item if $p\in e$, $L_p$ is a graph with two vertices and
four edges joining one vertex to the other, so $\chi(L_p)= 2 - 4 =
-2$;
\item if $p\in d$, $L_p$ is two points, so $\chi(L_p)= 2$;
\item if $p = (0,0,0)$, $L_p$ is a graph with two vertices and two edges
which are
both loops at the same vertex ($L_p$ is the union of a figure eight and a
point), so
$\chi(L_p)= 2 - 2 = 0$.
\end{itemize}

\begin{figure}[h]
\centering
\includegraphics[scale=.5]{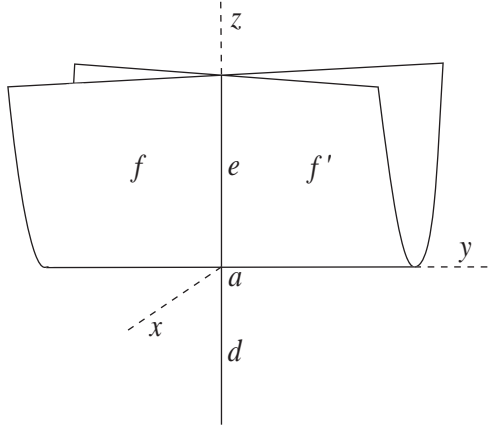}
\caption{The Cartan umbrella}
\end{figure}

On the other hand, the {\em Whitney umbrella} $X=\{(x,y,z)\ |\ x^2 = zy^2,\
z\geq
0\}$ is not algebraic, because the link of the origin is a figure eight,
which has
Euler characteristic $-1$. (The Whitney umbrella is the image of the
algebraic map
$(u,v)\mapsto (uv,u,v^2)$.)

Now $\chi(L_p)$ is a topological invariant, since it can be computed from
the local
homology of $X$ at $p$:
$$\chi(L_p)= 1 - \sum_k (-1)^k \operatorname{rank} H_k(X,X-\{p\}).$$
Thus the Whitney umbrella is not  homeomorphic to a real algebraic set.

The known higher topological obstructions for a semialgebraic set to be
algebraic are
similar to Sullivan's obstruction. They are local---in fact they are mod 2
invariants of the links $L_p$---and they are combinatorial in the sense
that they
depend only on the dimensions and incidence  of the cells of a
semialgebraic cell structure on $X$.

Parusi\'nski and the author \cite{mccrory-parusinski1} have found a way to describe such
obstructions using the ring
of constructible functions on $X$ and an operator on this ring generalizing
the Euler
characteristic of the link. Here is a brief
description of our method.

Let $X$ be a semialgebraic set. An
integer-valued function
$\varphi:X\to \Z$ is {\em constructible} if
it is a finite combination of characteristic functions of semialgebraic subsets
$X_i$ of $X$ with integer coefficients: $\varphi = \sum_i n_i\mathbf
1_{X_i}$. It
follows that
there exists a semialgebraic cell
complex $\mathcal K$ on $X$ such that $\varphi$ is constant on each
cell of
$\mathcal K$. (A semialgebraic cell complex is a locally topologically trivial
semialgebraic stratification such that each stratum is
homeomorphic to a Euclidean space and has compact closure.)

The set of constructible functions on $X$ is a ring, with $(\varphi +\psi)(p)
= \varphi(p) +\psi(p)$ and $(\varphi \times\psi)(p) = \varphi(p)
\times\psi(p)$. If
$\varphi = \sum_i n_i\mathbf 1_{X_i}$ with the sets $X_i$ compact,  the
{\em Euler
integral} of $\varphi$ is defined by
$\int_X\varphi =  \sum_i n_i \chi(X_i)$, where $\chi$ is the Euler
characteristic.
If
$\varphi$ is constant on each cell of $\mathcal K$ and $\varphi$ has compact
support, then
$$\int_X\varphi = \sum_{C\in\mathcal K} (-1)^{\dim C} \varphi(C),$$
where $\varphi(C)$ is the value of $\varphi$ on the cell $C$.

The {\em link operator} $\link$ on the ring of constructible functions on
$X$ is
the Euler integral of the restriction of $\varphi$ to the link of $X$ at $p$,
$$(\link\varphi)(p) = \int_{L_p}\varphi.$$

Sullivan's theorem says that if $X$ is an algebraic set with characteristic function
$\mathbf 1_X$, the values of the function
$\link
\mathbf 1_X$ are even. In other words, if $\hlink = \frac 12 \link$ then
$\hlink\mathbf 1_X$ is integer-valued (and hence $\hlink\mathbf 1_X$ is
a new constructible function on $X$).

 The main result of
\cite{mccrory-parusinski1} (Theorem 2.5, p. 536) implies that if $X$ is an
algebraic
set, then {\em all} the functions obtained from $\mathbf 1_X$ using the
arithmetic
operations $+$, $-$,
$\times$, and the operator $\hlink$ are integer-valued. Furthermore this
property is
topological (\cite{mccrory-parusinski1} A.7, p. 550).

The following statement is a corollary: If $Y$ is an algebraic
link (the link of a point in an algebraic set), then all the
functions obtained from $\mathbf 1_Y$ using the arithmetic operations $+$, $-$,
$\times$, and the operator $\hlink$ are integer-valued and have even Euler
integral.

\begin{figure}[h]
\centering
\includegraphics[scale=.5]{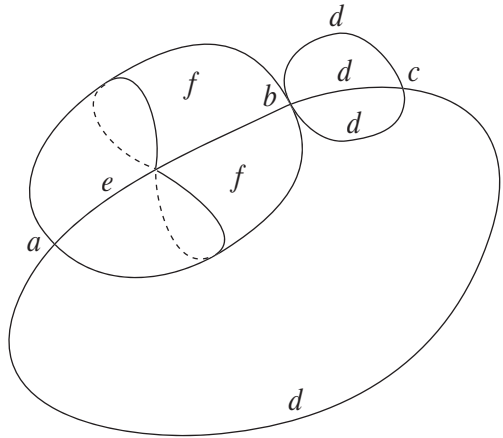}
\caption{The space Y}
\end{figure}

Now we apply our method to Akbulut and King's example. Let
$Y$ be the topological space shown in Figure 2.
The space $Y$ is homeomorphic to an algebraic set in projective 3-space,
the union
of the projective Cartan umbrella $wx^2 = zy^2$ and the plane conic
$(w+z)^2 + x^2
= z^2$, $y=0$. Realized in this way, $Y$ has a semialgebraic cell structure
with
three 0-cells, five 1-cells, and two 2-cells. The three 0-cells are
labelled $a,
b, c$ in the figure. Four of the 1-cells are labelled $d$, and one of the
1-cells
(incident to the 2-cells) is labelled $e$. The two 2-cells are labelled $f$.

The table below summarizes the computation of a constructible function $\psi$
obtained from $\mathbf 1_Y$ using the arithmetic operations $+, -, \times$,
and the
operator $\hlink$, such that $\int_Y \psi$ is odd.

The strategy for finding a suitable function $\psi$ is based on the following
observations. We define an operator $\ahlink$ by $\ahlink\xi = \xi -
\hlink\xi$. The interesting fact $\hlink\hlink =
\hlink$ implies that $\ahlink\ahlink =
\ahlink$, $\hlink\ahlink =
0$, and $\ahlink\hlink =
0$.
If $\xi$ has even-dimensional support then the support of $\hlink\xi$ is
of lower dimension. If If $\xi$ has odd-dimensional support then the support of
$\ahlink\xi$ is of lower dimension.

Let $d$ be the dimension of $Y$. For decreasing dimensions $k =
d,d-1,d-2,\dots,$ we
search through the functions with support dimension $k$ which are obtained from
$\mathbf 1_Y$ using $+, -, \times, \hlink$, until we find one with odd Euler
integral.

In the example at hand, since the dimension of $Y$ is 2 we start with
$\varphi =
\hlink \mathbf 1_Y$, which has support dimension 1, and we apply
$\ahlink$ to the powers of $\varphi$ to get functions with support dimension
zero.

Computation of the obstruction:

\smallskip

\begin{center}
\begin{tabular}{| c |  c | c | c | c | c | c || c |}
\hline
 & $a$ & $b$ & $c$ & $d$ & $e$ & $f$ & $\int$\\\hline
$\mathbf 1_Y$ & 1 & 1 & 1 & 1 & 1 & 1 & 0 \\
$\varphi = \hlink\mathbf 1_Y$ & 0 & 1 & 2 & 1 & $-1$ & 0 & 0 \\
$\varphi^2$ & 0 & 1 & 4 & 1 & 1 & 0 & 0 \\
$\hlink\varphi^2$ & 1 & 2 & 2 & 1 & 1 & 0 & 0 \\
$\ahlink\varphi^2= \varphi^2 - \hlink\varphi^2$ & $-1$ & $-1$ & 2 & 0 &
0 & 0 & 0 \\
$\psi =\varphi\ahlink\varphi^2$ & 0 & $-1$ & 4 & 0 & 0 & 0 & 3 \\\hline
\end{tabular}
\end{center}

\medskip

In the table the columns $a,b,c,d,e,f$ correspond to the types of cells of
$Y$. Two cells have the same label if they have the same local topology at an
interior point. This implies that any function obtained from $\mathbf
1_Y$ using $+, -, \times, \hlink$ will be constant on all cells with the same
label. The rows of the table are such functions, and the entries of the
table give
the values of the function on each type of cell.

The last column is the Euler
integral of the function. To compute the Euler integral we take into
account that
there are four cells of type $d$ and two cells of type $f$. Thus, for
example, for
the first row we have $\int \mathbf 1_Y = (1+1+1)-(4\cdot 1 + 1) + (2\cdot
1) = 0$,
and for the second row we have $\int \varphi = (0+1+2)-(4\cdot 1 + (-1)) +
(2\cdot
0) = 0$.

The computation of $\varphi = \hlink\mathbf 1_Y$ is similar to our first
example, the Cartan umbrella. The link of the point $a$ is the union of a
figure eight and a point, which has Euler characteristic 0, so $\varphi(a)
= 0$. The
link of $b$ is the union of a figure eight and three points, which has Euler
characteristic 2, so $\varphi(b) = 1$. The link of $c$ is four points so
$\varphi(c) = 2$. The link of a point of type $d$ is two points so
$\varphi(d) = 1$.
The link of a point of type $e$ is a graph with two vertices and four edges
from
one vertex to the other, which has Euler characteristic $-2$, so
$\varphi(e) = -1$.
And the link of a point of type $f$ is a circle, which has Euler
characteristic 0, so
$\varphi(f) = 0$.

The computation of $\hlink\varphi^2$ is similar. For example, the link of $a$
is the union of a figure eight and a point. The point is of type $d$, at which
$\varphi^2$ takes the value 1. The vertex of the figure eight is of type
$e$, at
which
$\varphi^2$ also takes the value 1. The 1-cells of the figure eight are of
type $f$,
at which $\varphi^2$ takes the value 0. Thus the Euler integral of
$\varphi^2$ on
the link of $a$ is $(1+1) - 0 = 2$, so $\hlink\varphi^2(a) = 1$.

The result of our computation is
$$\int_Y\psi = 3,\ \psi = (\hlink \mathbf 1_Y)\big((\hlink \mathbf
1_Y)^2-\hlink(\hlink \mathbf 1_Y)^2\big).$$
Thus $Y$ is not homeomorphic to an algebraic link.

Akbulut and King's example $X$ is the {\em suspension} of $Y$, the
union of
two cones with base $Y$. Now
$X$ can be realized as a semialgebraic set in a Euclidean space, with $Y$
the link
of the vertex of either cone. Since the Euler
characteristic of $Y$ is 0, $X$ satisfies Sullivan's condition. Since $Y$
is not
homeomorphic to an algebraic link, it follows that $X$ is not homeomorphic
to an
algebraic set.

A more direct way to prove that $X$ is not homeomorphic to an algebraic set
is to
find a function $\eta$ obtained from $\mathbf 1_X$ using
$+, -, \times, \hlink$ such that $\eta$ is not integer-valued. It is not hard
to see that the operators $\hlink$ and $\ahlink$ are interchanged by
restriction to a link (see \cite{mccrory-parusinski2} 1.3(d), p. 499). More
precisely, if $Y$ is the link of $p$ in $X$, and
$\xi$ is a constructible function on $X$, then
$$(\hlink\xi)|_Y = \ahlink(\xi|_Y),\
(\ahlink\xi)|_Y = \hlink(\xi|_Y).$$
We apply this to the Akbulut-King example. Since $(\mathbf 1_X)|_Y =
\mathbf 1_Y$, we
have $(\ahlink\mathbf 1_X)|_Y = \hlink\mathbf 1_Y$. In other words, if
$\varphi' =
\ahlink\mathbf 1_X$ then $\varphi'|_Y = \varphi$. Thus
$\varphi'\hlink(\varphi')^2|_Y
= \varphi\ahlink\varphi^2$. So if $p$ is one of the two suspension vertices
of $X$
then $\hlink(\varphi'\hlink(\varphi')^2)(p) = \frac 12\int_Y
\varphi\ahlink\varphi^2
=\frac 32$. The function $\eta = \ahlink(\varphi'\hlink(\varphi')^2)$ is
supported
at the two vertices of $X$, and at each of these vertices $\eta(p) = -\frac
32$.

Akbulut and King \cite{akbulut-king} analyze five local topological
obstructions for
a 3-dimensional semialgebraic set to be algebraic---the local mod 2 Euler
characteristic and four new invariants. They prove that these five
obstructions are
independent, and that the vanishing of these obstructions is a sufficient
condition
for a compact 3-dimensional semialgebraic set to be homeomorphic to an
algebraic set.

The situation in dimension four is not so simple. There are at least
$2^{43}-43$
independent local topological obstructions \cite{mccrory-parusinski2}, but
it is not
known whether the vanishing of all these obstructions is a sufficient
condition for a
compact 4-dimensional semialgebraic set to be homeomorphic to an algebraic set.

Significant
progress toward the description of a set of sufficient conditions in dimension
four has been made by Michelle Previte \cite{previte} using Akbulut and King's
{\em resolution towers}. However, the relation between resolution tower
obstructions
and our constructible function invariants remains a mystery.

\end{document}